\documentclass{article}
\usepackage{enumerate}
\usepackage[toc,page]{appendix}
\usepackage{amsmath,amssymb,amsthm}
\usepackage{graphicx}

\usepackage[dvipsnames]{xcolor}
\usepackage{hyperref}
\hypersetup{
    colorlinks=true,
    linkcolor=cyan!80!black,
    citecolor=MidnightBlue,
    urlcolor=magenta,
}

\theoremstyle{plain}
\newtheorem{theorem}{Theorem}

\newtheorem{lemma}{Lemma}

\theoremstyle{definition}

\newtheorem{definition}{Definition}

\let\tilde\widetilde

\title{Overlap distribution of spherical spin glass models with general eigenvalue distribution of the interaction matrix}
\author{Debapratim Banerjee and Debabrata Jana\\
Department of Mathematics \\
    Ashoka University, India\\
 debapratim.banerjee@ashoka.edu.in \\
 debabrata.jana05@gmail.com }
\date{September 2025}

\begin{document}

\maketitle
\begin{abstract}
    In this paper, we show that the replica symmetry of the Gibbs measure of spherical spin systems is a property of the eigenvalue spacing at the edge of the interaction matrix. In particular, our interaction matrix has \textbf{two} large outlier eigenvalues with mutual distance $\frac{c}{n}$. The empirical measure of the rest of the eigenvalues is close to the semicircular law with some rigidity conditions. We prove that in this scenario the overlap distribution of two independent samples from the Gibbs measure has a continuous density at a low enough temperature. Hence, the model is a full replica symmetry-breaking model. One might compare this result with only one outlier eigenvalue. This model comes for the Sherrington-Kirkpatrick model with Curie-Weiss interaction in the ferromagnetic case. Here, it is well known that the model is replica symmetric, although the free energy limit of this model is the same as the free energy limit of our model. In our limited understanding, we believe that this kind of phenomenon cannot be explained by the Parisi approach.  
\end{abstract}
\section{Introduction}
The study of the overlap distribution of i.i.d. samples from the Gibbs measure of spin systems is of fundamental interest in the literature. In fact, following the amazing idea of Parisi, the literature of Sherrington Kirkpatrick model based on the Parisi approach is about relating the overlap distribution to the free energy of the system. In the last few years, we have seen remarkable breakthroughs in the literature of spin glasses. This started with the breakthrough paper of Talagrand \cite{Tal05} where he proved that the free energy limit of the Sherrington Kirkpatrick model at any temperature is given by the Parisi formula. On the basis of this paper, a large amount of work was done on the Sherrington Kirkpatrick model. We shall not be able to discuss this in this introduction. The Parisi formula says that the free energy of the Sherrington Kirkpatrick model is given by an optimization procedure where the objective function is a functional of distribution functions from $[0,1] \to [0,1]$. It turns out that the optimizer is unique and is called the Parisi measure \cite{Tal05}. The Parisi measure is closely related to the limiting distribution of overlaps of two independent samples under the Gibbs measure. See \cite{panchenko2007overlap} for the spherical case. To our limited knowledge, the conjecture that the Ising $2$ spin Sherrington Kirkpatrick model is full RSB at low enough temperature is still widely open. However, the spherical $2$ spin case is fully solved. It turns out that the spherical $2$ spin SK model is replica symmetric at any temperature. See \cite{panchenko2007overlap} page 2323 for a reference. There have been intensive studies on spherical mixed spin case as well. One might look at \cite{auffinger2019existence} and Proposition 1 of \cite{subag2021following} and the references there in. However, in our limited understanding, the spherical $2$ spin case with general interaction matrix is not yet studied in the literature. In this paper, we initiate the study. In this paper, we construct an interaction matrix where the spherical $2$ spin model becomes full RSB at low enough but non zero temperature. This result is surprising in the following two senses. First of all, the limiting free energy is no longer a function of the overlap distribution. In particular, there are infinitely many candidates for overlap distributions which will produce the same free energy and limiting overlap distribution turns out to be one of them. Secondly, the model again becomes replica symmetric at $0$ temperature provided $c>0$. Here $\frac{c}{n}$ is the difference between the two largest eigenvalues. The proof of this paper originated from \cite{Ban25} where the first author of this paper introduced a novel approach to study the free energy of spherical $2$ spin models directly. We hope the observation of this paper will have some importance in the study of $2$ spin spherical spin glass models.     

\section{The model description}
In this section, we give the formal definition of the model.
Firstly, we define the uniform probability distribution on the unit sphere. 
\begin{definition}\label{def:uni}(Uniform distribution on the sphere)
For any $n \ge 2$, define $\mathbb{S}_{n-1}$ as the set of all $\underline{\sigma}=(\sigma_{1},\ldots, \sigma_{n})' \in \mathbb{R}^{n}$ such that $||\underline{\sigma}||^{2}=1$. $\mathbb{S}_{n-1}$ is called the unit sphere in $\mathbb{R}^{n}$. Let $\mu_{n}$ be the uniform probability measure on $\mathbb{S}_{n-1}$.
\end{definition}
\noindent 
Next, we define the Hamiltonian of the Gibbs measure of the spin systems. 
\begin{definition}\label{def:ham}(Hamiltonian)
    Let $X_{n}$ be a matrix of dimension $n \times n$. For any $\underline{\sigma} \in \mathbb{S}_{n-1}$, the Hamiltonian of the model is defined as 
    \[
    H_{n}(\underline{\sigma})= \frac{n}{2} \underline{\sigma}' X_{n} \underline{\sigma}.
    \]
\end{definition}
\begin{definition}(Partition function and free energy)\label{def:partitionfree}
    The partition function of the model at an inverse temperature $\beta >0$, denoted as $Z_{n}(\beta)$, is defined as follows:
    \[
      Z_{n}(\beta)= \int_{\mathbb{S}_{n-1}} \exp\left\{ \beta H_{n}(\underline{\sigma}) \right\} d\mu_{n}(\underline{\sigma}).
    \]
    The Free energy, $F_{n}(\beta)$ is defined as $\log (Z_{n}(\beta))$.
\end{definition}
\begin{definition}\label{def:gibbs}(Gibbs measure)
  Finally, the Gibbs measure $\mu_{n,\beta}^{\mathrm{Gibbs}}$ is a probability measure on $\mathbb{S}_{n-1}$ is defined as follows: 
  \[
   d\mu_{n,\beta}^{\mathrm{Gibbs}}(\underline{\sigma})= \frac{\exp\left\{ \beta H_{n}(\underline{\sigma}) \right\}}{Z_{n}(\beta)} d\mu_{n}(\underline{\sigma}).
  \]
\end{definition}
In this paper, we shall be interested in the measure $\mu_{n,\beta}^{\mathrm{Gibbs}}$ when $n \to \infty$ and $\beta$ is a fixed but large number.
\begin{definition}\label{def:overlap}(Overlap distribution)
    Let $\underline{\sigma}^{(1)}$ and $\underline{\sigma}^{(2)}$ be two samples from the Gibbs distribution $\mu_{n,\beta}^{\mathrm{Gibbs}}$. The overlap between $\underline{\sigma}^{(1)}$ and $\underline{\sigma}^{(2)}$ shall be denoted by $\mathrm{Ov}(\underline{\sigma}^{(1)},\underline{\sigma}^{(2)})$ and defined as 
    \[
     \mathrm{Ov}(\underline{\sigma}^{(1)},\underline{\sigma}^{(2)}) = \langle \underline{\sigma}^{(1)},\underline{\sigma}^{(2)} \rangle = \sum_{i=1}^{n} \sigma^{(1)}_{i} \sigma^{(2)}_{i}.
    \]
    The measure corresponding to the distribution of $\mathrm{Ov}(\underline{\sigma}^{(1)},\underline{\sigma}^{(2)})$ will be of central interest of the current paper.
\end{definition}
\begin{definition}\label{def:replica}(Replica symmetry)
    Observe that for our model, $\mathrm{Ov}(\underline{\sigma}^{(1)},\underline{\sigma}^{(2)})$ is symmetric around $0$. If $\left| \mathrm{Ov}(\underline{\sigma}^{(1)},\underline{\sigma}^{(2)}) \right|$ is asymptotically (i.e. $n \to \infty$) concentrated around a single value, then the model is called replica symmetric. If it is asymptotically concentrated around $k+1$ values for $k \in \mathbb{N}$, then the model is called $k-$RSB (replica symmetry breaking). In all the other cases, the model is called full-RSB. 
\end{definition}
\noindent 
In the context of spin systems, an important model is the spherical Sherrington- Kirkpatrick model where $X_{n}$ is taken to be the G.O.E. matrix (Gaussian orthogonal ensemble). For the Sherrington- Kirkpatrick model with Curie Weiss interaction, one takes $X_{n}= W_{n}+ \frac{J}{n} \mathbf{11}'$ where $J >0$ and $\mathbf{1}$ is the $n \times 1$ vector of all ones. We shall formally define these matrices next.
\begin{definition}(GOE matrix)
    Let $X_{n}= \frac{1}{\sqrt{n}}\left(  x_{i,j} \right)_{1\le i \le j \le n}$ be a symmetric matrix of dimension $n \times n$. It is called a GOE matrix (Gaussian orthogonal ensemble) if $x_{i,j}\sim N(0,1)$ for $1\le i < j \le n$, $x_{i,i} \sim N(0,2)$ for $1 \le i \le n$ and the random variables are independent.
\end{definition}
Let $\lambda_{1}\ge \ldots \ge \lambda_{n}$ be the eigenvalues of the matrix $X_{n}$. It is a classical result in random matrix theory that 
\[
\frac{1}{n}\sum_{i=1} \delta_{\lambda_{i}}
\]
converges weakly to the semicircular distribution in the almost sure sense. We now formally define the semicircular law.
\begin{definition}\label{def:sc}
    The semicircular law is a probability distribution on $[-2,2]$ and it's probability density function is 
    \begin{equation}
        f(x)= \frac{1}{2\pi} \sqrt{4-x^2}~ \mathbb{I}_{|x|\le 2}.
    \end{equation}
    In this paper, we denote the probability measure of the semicircular law by $\mu_{\mathrm{sc}}$.
\end{definition}
\begin{definition}(Spiked random matrices)
    When the matrix $X_{n}$ can be written as $W_{n}+V$ where $V$ is a symmetric deterministic matrix of some fixed rank $r$ and $W_{n}$ is the G.O.E. matrix, $X_{n}$ is called a spiked matrix. We shall be interested in the following special case: $X_{n}= W_{n}+ \frac{J}{n} \mathbf{11}'$. 
\end{definition}

\noindent 
It is well known in the random matrix literature that when $J>1$, the largest eigenvalue of $X_{n}$ converges to $J+\frac{1}{J}$ and the empirical measure of $\lambda_{2}\ge \ldots\ge \lambda_{n}$ converges weakly almost surely to the semicircular law. 
\begin{definition}(Our interaction matrix)\label{def:ourint}
The matrices we shall consider in this paper will have eigenvalues of the following kind: Suppose $X_{n}$ is matrix with eigenvalues $\lambda_{1}\ge \ldots \ge \lambda_{n}$ where $\lambda_{1}= J+\frac{1}{J}$, $\lambda_{2}= \lambda_{1}- \frac{c}{n}$  and the empirical measure of $\lambda_{2}\ge \ldots\ge \lambda_{n}$ converges weakly almost surely to the semicircular law with some rigidity conditions.
\end{definition}
\section{Main results}
The main results of this paper are as follows:
\begin{theorem}\label{thm:free}
    Suppose the interaction matrix $X_{n}$ satisfies the properties in Definition \eqref{def:ourint}, $\beta J >1$ and $J>1$. Then 
    \begin{equation}
        \frac{F_{n}(\beta)}{n} \to \frac{\beta}{2}\left( J+ \frac{1}{J}  \right)\left( 1- \frac{1}{\beta J} \right)+ \frac{1}{2}\log \left( \frac{1}{\beta J} \right) + \frac{1}{4J^2} ~~\text{almost surely}.
    \end{equation}
    This result is same as the free energy limit obtained by \cite{baik2017fluctuations} in Theorem 1.4 (1.10) (iii).
\end{theorem}
\begin{theorem}\label{thm:over}
    Suppose the interaction matrix $X_{n}$ satisfies the properties in Definition \eqref{def:ourint}, $\beta J >1$ and $J>1$. Let $\underline{\sigma}^{(1)}$ and $\underline{\sigma}^{(2)}$ be two independent samples from the probability distribution $\mu_{n,\beta}^{\mathrm{Gibbs}}$. Then as $n \to \infty$, 
    \begin{equation}\label{eq:resultden}
     \mathrm{Ov}\left( \underline{\sigma}^{(1)}, \underline{\sigma}^{(2)} \right) \stackrel{d}{\to} \sqrt{s_{1}s_{2}} \zeta_{1} + \sqrt{(\hat{r}_{0}- s_{1})(\hat{r}_{0}-s_{2})} \zeta_{2} ~~ \text{almost surely,}
    \end{equation}
    with $\hat{r}_{0}=\left( 1- \frac{1}{\beta J} \right)$. Here $\zeta_{1}, \zeta_{2}$ are i.i.d. Rademacher random variables and $s_{1}$, $s_{2}$ are i.i.d. random variables with density:
    \[
     f_{s}(x) \propto \left\{
       \begin{array}{cc}
         \frac{1}{\sqrt{x(\hat{r}_{0}-x)}} \exp\{ -\frac{c\beta}{2} x \} & \text{whenever } 0 \le x \le \hat{r}_{0} \\
           0 & \text{otherwise.}
       \end{array}
     \right.
    \]
\end{theorem}
\section{Dirichlet distributions}
The proof of Theorem \eqref{thm:free} requires some properties of Dirichlet distributions. So, we first introduce them.  
\begin{definition}(Dirichlet distributions)
Dirichlet distributions are a family of multivariate distributions parametrized by a vector $\underline{\alpha}=(\alpha_{1},\ldots, \alpha_{k})$. Suppose $\underline{Y}=(Y_{1}, \ldots, Y_{k}) \sim \mathrm{Dir}(\alpha_{1},\ldots, \alpha_{k})$, then $Y$ is supported on the simplex $S_{k}:=\{0\le y_{i} \le 1\} \cap \{ \sum_{i=1}^{k} y_{i} =1\}$. The p.d.f. (probability density function) at $(y_{1},\ldots , y_{k})$ is given by 
\begin{equation}\label{eq:dirden}
    f_{\alpha_{1},\ldots, \alpha_{k}}(y_{1},\ldots, y_{k})= \frac{\Gamma(\alpha_{1}+ \ldots + \alpha_{k})}{\prod_{i=1}^{k} \Gamma(\alpha_{i})} \prod_{i=1}^{k} y_{i}^{\alpha_{i}-1}.
\end{equation}
\end{definition}
We would require the following property of the Dirichlet distribution: 
\begin{lemma}\label{lem:diruse}
Suppose $X_{1},\ldots, X_{k}$ are independent $\Gamma(\alpha_{1},\tau),\ldots, \Gamma(\alpha_{k},\tau)$ respectively, then $(Y_{1},\ldots, Y_{k}) \sim \mathrm{Dir}(\alpha_{1},\ldots, \alpha_{k})$, where $Y_{i}:=  \frac{X_{i}}{X_{1}+ \ldots + X_{k}}$.
\end{lemma}
The proof of Lemma \eqref{lem:diruse} is simple and uses a multivariate change of variable formula. So we omit it.
We know that $\chi^{2}(2\alpha)$ distribution is same as  $\Gamma(\alpha,2)$ distribution. Hence, in Lemma \eqref{lem:diruse}, we can take $X_{1},\ldots, X_{k}$ independent $\chi^{2}(2\alpha_{1}), \ldots \chi^{2}(2\alpha_{k})$.
\section{Proofs of Theorems \eqref{thm:free} and \eqref{thm:over}}
\begin{proof}{Proof of Theorem \eqref{thm:free}.}
The proof of the theorem $\eqref{thm:free}$ contains the following four parts: (i) Diagonalization,\,(ii) Discretization,\,(iii) Optimization, and\,(iv) Laplace integration method. We shall complete these step one by one.
\smallskip

\textbf{(i)\,Diagonalization:} This step is straightforward. As $\underline{\eta}$ is uniformly distributed over $\mathbb{S}_{n-1},\, H\underline{\eta}$ is also uniformly distributed over $\mathbb{S}_{n-1}$ for any orthogonal matrix $H$. As a consequence,
\begin{eqnarray}
     Z_{n}(\beta)&=& \int_{\mathbb{S}_{n-1}} \exp\left\{ \beta H_{n}(\underline{\eta}) \right\} d\mu_{n}(\underline{\eta})\nonumber\\
      &=& \int_{\mathbb{S}_{n-1}} \exp\left\{\frac{n\beta}{2}\underline{\eta^{\prime}} H_{n}^\prime\Lambda_n H_{n}\underline{\eta} \right\} d\mu_{n}(\underline{\eta})\\
      &=& \int_{\mathbb{S}_{n-1}} \exp\left\{ \frac{n\beta}{2}\sum_{i=1}^{n}\lambda_i\eta_i^2\right\} d\mu_{n}(\underline{\eta})\nonumber
\end{eqnarray}
 Here $X_n=H_{n}^\prime\Lambda_n H_{n}$ is spectral decomposition of the interaction matrix $X_n$ 
 and $\lambda_1 (=J+\frac{1}{J})\ge\lambda_{2} (= J+\frac{1}{J}- \frac{c}{n})\ge\lambda_3\ge \ldots\ge \lambda_{n}$ are eigen values of $X_n$.
 \smallskip

 \textbf{(ii) Discretization:}
 We take finite but 
 large enough $K$ and make a discrete approximation of the semicircular law. Let $\tilde{\lambda}_{1}\ge\tilde{\lambda}_{2}\ge \ldots\ge \tilde{\lambda}_{k+1}=-2$ be such that \[\int_{\tilde{\lambda}_{j+1}}^{\tilde{\lambda}_{j}}d\mu_{sc}(x)=\frac{1}{K}
\]
Fixing $\varepsilon>0,$ we take $K$ large enough so that
\begin{eqnarray}\label{eq:rigid}
\sup_{1\le i\le K}\sup_{\frac{(i-1)(n-2)}{K}+3\le j\le \frac{i(n-2)}{K}+2}\Big|\lambda_j-\tilde{\lambda}_{i}\Big|\le \varepsilon
\end{eqnarray}
This is ensured from the eigenvalue rigidity of the GOE matrix. For example see \cite{erdHos2012rigidity}. From \eqref{eq:rigid}, observe that
\begin{eqnarray}\label{}
\left| \frac{n\beta}{2}\sum_{j=3}^{n}\lambda_j\eta_j^2 - \frac{n\beta}{2}\sum_{i=1}^{n}\tilde{\lambda}_i\left(\sum_{j= \frac{(i-1)(n-2)}{K}+3}^{\frac{i(n-2)}{K}+2}\eta_j^2\right)\right|\le \frac{n\beta\varepsilon}{2}
\end{eqnarray}
Hence we shall work with 
\begin{eqnarray}\label{eq:dpart}
\int_{\mathbb{S}_{n-1}} \exp\left\{ \frac{n\beta}{2}\left((J+\frac{1}{J})(\eta_1^2+\eta_2^2)+\sum_{i=1}^{K}\tilde{\lambda}_i\left(\sum_{j=\frac{(i-1)(n-2)}{K}+3}^{\frac{i(n-2)}{K}+2}\eta_j^2\right) \right)\right\} d\mu_{n}(\underline{\eta})
\end{eqnarray}
instead of actual partition function. Recall that $(\eta_1,\ldots,\eta_n)^\prime$ has the same distribution as $(\frac{Y_1}{||Y_1||},\ldots,\frac{Y_n}{||Y_n||})^\prime$ where $Y=(Y_1,\ldots,Y_n)^\prime \sim N_n(0,I_n)$. As a consequence the vector $\left(\eta_1^2+\eta_2^2,\, \sum_{j=\frac{(i-1)(n-2)}{K}+3}^{\frac{i(n-2)}{K}+2}\eta_j^2\right)_{1\le i\le k} \sim \text{Dir}(1, \frac{n-2}{2k},\ldots,\frac{n-2}{2k})$. Hence, \eqref{eq:dpart} can be written as 
\[\text{E}\left[\text{exp}\left\{ \frac{n\beta}{2}\sum_{i=0}^{K}\tilde{\lambda}_ir_i\right\}\right]\]
where $(r_0,\ldots ,r_K)\sim \text{Dir}(1, \frac{n-2}{2k},\ldots,\frac{n-2}{2k})$. Hence, 
\begin{eqnarray}\label{eq:expt}
\text{E}\left[\text{exp}\left\{ \frac{n\beta}{2}\sum_{i=0}^{K}\tilde{\lambda}_ir_i\right\}\right]=\int_{S_K}\text{exp}\left\{ \frac{n\beta}{2}\sum_{i=0}^{K}\tilde{\lambda}_ir_i\right\}f_{1, \frac{n-2}{2k},\ldots,\frac{n-2}{2k}}(r_0,\ldots, r_K)\,dr_0\ldots dr_K
\end{eqnarray}
This concludes the discretization.
\smallskip

\textbf{(iii) Optimization:} We now come to the third step. First of all we do some asymptotic on the p.d.f of $\text{Dir}(1, \frac{n-2}{2k},\ldots,\frac{n-2}{2k})$. From \eqref{eq:dirden} we know that 
\begin{eqnarray}\label{}
f_{1, \frac{n-2}{2k},\ldots,\frac{n-2}{2k}}(r_0,\ldots, r_K)=\frac{\Gamma(\frac{n}{2})}{ \Gamma(\frac{n-2}{2K})^K} \prod_{i=1}^{k} r_{i}^{\frac{n-2}{2K}-1}
\end{eqnarray}
We approximate $\Gamma(\frac{n}{2})$ by $(\frac{n}{2})^{\frac{n}{2}}\text{exp}\{-\frac{n}{2}+o(n)\}$ and $\Gamma(\frac{n-2}{2K})$ by $(\frac{n-2}{2K})^{\frac{n-2}{2K}}\text{exp}\{-\frac{n-2}{2K}+o(n)\}.$
These are estimate the obtained by Stirling's approximation. Hence,
\begin{eqnarray}\label{}
\frac{\Gamma(\frac{n}{2})}{ \Gamma(\frac{n-2}{2K})^K}=(K)^{\frac{n-2}{2}}\text{exp}\{o(n)\}=\text{exp}\left\{\frac{n-2}{2}\text{log}K+o(n)\right\}
\end{eqnarray}
We shall now consider the integrand in the r.h.s of \eqref{eq:expt} and shall optimize with respect to $r_0,\ldots,r_K$. The integrand in the r.h.s of \eqref{eq:expt} becomes:
\begin{eqnarray}\label{eq:maxim}
\text{exp}\left\{\frac{n\beta}{2}\sum_{i=0}^{K}\tilde{\lambda}_ir_i+\frac{n}{2}\text{log}K+\left(\frac{n-2}{2K}-1\right)\sum_{i=1}^K\text{log}r_i\right\}
\end{eqnarray}
The rest of this step is devoted to maximize \eqref{eq:maxim} where $r_0,\ldots,r_K$ varies over $S_K$.
Hence it is enough to find 
\begin{eqnarray}\label{eq:maxim1}
\max_{(r_0,\ldots,r_K)\in S_K} F(r_{0},\ldots, r_{K}):= \frac{n\beta}{2}\sum_{i=0}^{K}\tilde{\lambda}_ir_i+\frac{n}{2}\text{log}K+\frac{n}{2K}\sum_{i=1}^K\text{log}r_i
\end{eqnarray}
\eqref{eq:maxim1} is concave in $r_0,\ldots,r_k$. Hence the Lagrange multiplier method gives the maximum value when $r_0,\ldots,r_k$ varies over $S_K$. Consider the following function: 
\begin{eqnarray}\label{eq:lagran}
M(r_0,\ldots,r_k):=\beta\sum_{i=0}^{K}\tilde{\lambda}_ir_i+\text{log}K+\frac{1}{K}\sum_{i=1}^K\text{log}r_i-\Gamma\left(\sum_{i=0}^Kr_i-1\right)
\end{eqnarray}
Taking $\frac{\partial M}{\partial r_i}$ and putting it $0,$ we have 
\begin{eqnarray}\label{eq:partial}
\frac{\partial M}{\partial r_0}&=&\beta\tilde{\lambda}_0-\Gamma=0, \, \frac{\partial M}{\partial r_i}=\beta\tilde{\lambda}_0+\frac{1}{r_iK}-\Gamma=0,\,1\le i\le K \nonumber\\ 
\Gamma&=& \beta (J+\frac{1}{J}), \,\,\qquad r_i=\frac{1}{K(\Gamma-\beta\tilde{\lambda}_i)},\,1\le i\le K
\end{eqnarray}
The fact $\sum_{i=0}^Kr_{i}=1$ gives
\begin{eqnarray}\label{eq:conv}
\sum_{i=1}^{K}\frac{1}{K(\beta (J+\frac{1}{J})-\beta\tilde{\lambda}_i)}=1-r_0,
\end{eqnarray}
where $(J+\frac{1}{J})>\tilde{\lambda}_i,\, 1\le i\le K$. Observe that as $K$ becomes larger and larger, the l.h.s of \eqref{eq:conv} converges to $ -\frac{1}{\beta}S_{\mu_{sc}}(J+\frac{1}{J})$. We know that $-S_{\mu_{sc}}(J+\frac{1}{J})=\frac{1}{J}<1$. Thus, from \eqref{eq:conv} we get $1-r_0=\frac{1}{\beta J}+o_K(1)$ and $r_i=\frac{1}{K(\beta (J+\frac{1}{J})-\beta\tilde{\lambda}_i)}+o_K(1)$. This are the optimizing value of $r_i,\,0\le i\le K.$
Put $v_i=\frac{r_i}{1-r_0},\,1\le i\le K$. Then $\sum_{i=1}^Kv_i=1$. Hence at the optimizer, the value  of \eqref{eq:maxim1} becomes:
\begin{eqnarray}\label{eq:optmv}
\frac{n}{2}\left[\beta(J+\frac{1}{J})(1-\frac{1}{\beta J})+\sum_{i=1}^{K}\frac{1}{J}\tilde{\lambda}_iv_i+\text{log}K+\frac{1}{K}\sum_{i=1}^K\text{log}\frac{v_i}{\beta J}\right]+o(n) 
\end{eqnarray}
From the proof of high temperature case $(\frac{1}{J}<1)$ in \cite{Ban25},  Theorem $1$ (case-I), we get
\[\frac{n}{2}\left[\sum_{i=1}^{K}\frac{1}{J}\tilde{\lambda}_iv_i+\text{log}K+\frac{1}{K}\sum_{i=1}^K\text{log}v_i\right]=\frac{n}{4J^2}+o(n)\]
Putting above in  \ref{eq:optmv}, we have \ref{eq:optmv} is
\begin{eqnarray}\label{eq: optvalue}
\frac{n}{2}\beta(J+\frac{1}{J})(1-\frac{1}{\beta J})+\frac{1}{4J^2}+\frac{n}{2}\text{log}\frac{1}{\beta J}+o(n) 
\end{eqnarray}
This concludes the optimization step. We call the optimized value of \eqref{eq:maxim1} at
any $\beta$, $F_{opt}(\beta)$.
\smallskip

\textbf{(iv) Laplace's method of integration:}
This step is standard. However, for
completeness, we sketch it here. Firstly, the r.h.s. of \eqref{eq:expt} is clearly bounded above by
$\text{exp}\{F_{opt}(\beta)\}$. On the other hand, the quantity \eqref{eq:maxim1} is continuous in it’s arguments and the solution is in the interior of $S_K$. So for any $\varepsilon>0$, we can find
an open ball $B_{\varepsilon}$ around the optimizer so that the value of \eqref{eq:maxim1} is uniformly bigger
than $F_{opt}-\varepsilon$ on $B_{\varepsilon}$. So r.h.s. of \eqref{eq:expt} is bigger than $\text{exp}\{F_{opt}-n\varepsilon+ \text{log}(\text{Vol}(B_\varepsilon))\}$.
As $\text{log}(\text{Vol}(B_\varepsilon))$ is a finite quantity for every fixed $\varepsilon$, we get the result.
\end{proof}

In order to prove Theorem \eqref{thm:over}, we need the following technical result.
\begin{lemma}\label{lem:open}
Let $(\varepsilon , K)$ be as in the proof of Theorem \ref{thm:free}. Take $\tilde{\varepsilon} = \max\{ \varepsilon^{\frac{1}{1000}}, \left(\frac{1}{K}\right)^{\frac{1}{1000}} \}$. Let $(\hat{r}_{0}, \ldots, \hat{r}_{K})$ be the optimizer in Theorem \ref{thm:free}. Consider the set $S(\hat{r}_{0},\ldots , \hat{r}_{K})= \{(r_{0},\ldots, r_{K}) \in S_{K} ~|~ \max_{0 \le i \le K} |\hat{r}_{i}-r_{i}| \le \tilde{\varepsilon} \}$. Then on $S_{K} \backslash S(\hat{r}_{0},\ldots , \hat{r}_{K})$, uniformly $F(r_{0},\ldots, r_{K})\le F(\hat{r}_{0},\ldots, \hat{r}_{K})- cn\tilde{\varepsilon}^2$ for some constant $c>0$.  
\end{lemma}

With Lemma \eqref{lem:open} in hand, we now prove Theorem \eqref{thm:over}. 
\begin{proof}{Proof of Theorem \eqref{thm:over}}
Assuming Lemma \ref{lem:open}, we now prove Theorem \ref{thm:over}. The proof is based on a conditional argument. We condition on the value of $(r_{0}^{(j)},(\eta_{3}^{2})^{(j)},\ldots, (\eta_{n}^{2})^{(j)})$ for $j \in \{ 1,2 \}$. Observe that $(r_{0},\ldots, r_{K})$ is a function of $(r_{0},\eta_{3}^{2},\ldots, \eta_{n}^{2})$. In particular, 
\[
r_{i} = \sum_{j=\frac{(i-1)(n-2)}{K}+3}^{\frac{i(n-2)}{K}+2} \eta_{j}^{2}.
\]
Let $A=\{ \underline{\sigma} ~|~ (r_{0},\ldots, r_{K}) \in S(\hat{r}_{0}, \ldots, \hat{r}_{K})  \}$. Now, 
\begin{equation}
\begin{split}
&\mu_{n,\beta}^{\mathrm{Gibbs}}\left[ \mathrm{Ov}\left( \underline{\sigma}^{(1)}, \underline{\sigma}^{(2)} \right)  \in B\right]\\
&= \mu_{n,\beta}^{\mathrm{Gibbs}}\left[ \mathrm{Ov}\left( \underline{\sigma}^{(1)}, \underline{\sigma}^{(2)} \right)  \in B \cap \underline{\sigma}^{(1)} \in A \cap \underline{\sigma}^{(2)} \in A \right]\\
&~~~~~~+  \mu_{n,\beta}^{\mathrm{Gibbs}}\left[ \mathrm{Ov}\left( \underline{\sigma}^{(1)}, \underline{\sigma}^{(2)} \right)  \in B \cap\left( \underline{\sigma}^{(1)} \in A^{c} \cup \underline{\sigma}^{(2)} \in A^{c}\right) \right].
\end{split}
\end{equation}
Notice that 
\[
\mu_{n,\beta}^{\mathrm{Gibbs}}\left[ \mathrm{Ov}\left( \underline{\sigma}^{(1)}, \underline{\sigma}^{(2)} \right)  \in B \cap\left( \underline{\sigma}^{(1)} \in A^{c} \cup \underline{\sigma}^{(2)} \in A^{c}\right) \right] \le \mu_{n,\beta}^{\mathrm{Gibbs}}\left[ \left( \underline{\sigma}^{(1)} \in A^{c} \cup \underline{\sigma}^{(2)} \in A^{c}\right)\right]
\]
So we consider two cases: firstly, when $\underline{\sigma}^{(1)}$ and $\underline{\sigma}^{(2)}$ belong to $A$. Secondly, at least one of them is outside $A$. 
The second probability is bounded above by 
$
2\mu_{n,\beta}^{\mathrm{Gibbs}} \left[ \underline{\sigma}^{(1)} \in A^{c} \right].
$
Now 
\begin{equation}
\begin{split}
&\mu_{n,\beta}^{\mathrm{Gibbs}}\left[ \underline{\sigma} \in A^{c} \right]= \frac{1}{Z_{n}(\beta)}\mathrm{E}\left[ \exp\left\{ \beta  H_{n}(\underline{\sigma})\right\} \mathbb{I}_{\underline{\sigma} \in A^{c}} \right]\\
&\le \frac{\exp\{n\varepsilon\}}{\exp\{ - n\varepsilon\}}\frac{\exp\{n \sup_{(r_{0},\ldots, r_{K}) \in S_{K} \backslash S(\hat{r}_{0},\ldots, \hat{r}_{K}) } F(r_{0},\ldots, r_{K}) \}}{\exp \left\{ n F(\hat{r}_{0},\ldots, \hat{r}_{K}) \right\}}\\
&\le \exp\left\{ 2n\varepsilon - n \gamma \tilde{\varepsilon}^{2} \right\}\le \exp\left\{ -n \gamma' \tilde{\varepsilon}^{2}\right\}
\end{split}
\end{equation}
as $\tilde{\varepsilon}^2 \gg \varepsilon$.
Hence, we can ignore the second case. From now on, we shall assume that both $\underline{\sigma}^{(1)}$ and $\underline{\sigma}^{(2)}$ belong to $A$. This gives us $\left(r_{0}^{(j)},\ldots, r_{K}^{(j)}\right) \in S(\hat{r}_{0},\ldots, \hat{r}_{K})$. Hence,
$0 \le r_{i}^{(j)} \le \hat{r}_{i} + \tilde{\varepsilon} \le \tau \tilde{\varepsilon} $  for $1 \le i \le K$ implying  $r_{i}^{(j)}$'s are very small whenever $1 \le i \le K$ and $ \left| r_{0}^{(j)} - \hat{r}_{0}\right|\le \tilde{\varepsilon}$. 
We take a typical value of $(r_{0}^{(j)},(\eta_{3}^{2})^{(j)},\ldots, (\eta_{n}^{2})^{(j)})$ satisfying these properties. Firstly, we find the conditional density of $((\eta_{2}^{2})^{(j)},(\eta_{3}^{2})^{(j)},\ldots, (\eta_{n}^{2})^{(j)})$ given $(r_{0}^{(j)},(\eta_{3}^{2})^{(j)},\ldots, (\eta_{n}^{2})^{(j)})$ under the Gibbs measure. Since $((\eta_{3}^{2})^{(j)},\ldots, (\eta_{n}^{2})^{(j)})$ are already fixed, they remain the same. So, essentially, we need to find the conditional density of $\left( (\eta_{2}^{2})^{(j)} \left| (r_{0}^{(j)},(\eta_{3}^{2})^{(j)},\ldots, (\eta_{n}^{2})^{(j)}) \right. \right)$ under the Gibbs measure. A change of variable ($\left(r_{0}^{(j)}, (\eta_{2}^{2})^{(j)}\right) \mapsto \left( r_{0}^{(j)}- (\eta_{2}^{2})^{(j)}, (\eta_{2}^{2})^{(j)}  \right) $) will show that the conditional density of $\left((\eta_{2}^{2})^{(j)} \left| (r_{0}^{(j)},(\eta_{3}^{2})^{(j)},\ldots, (\eta_{n}^{2})^{(j)}) \right. \right)$ under the Gibbs measure is same as  
\begin{equation}\label{eq:densityfinale}
\begin{split}
&\frac{\exp\left\{ \frac{n\beta}{2} \left( J + \frac{1}{J} \right) r_{0}^{(j)}- \frac{c\beta}{2} (\eta_{2}^{2})^{(j)} + \frac{n\beta}{2} \sum_{i=3}^{n} \lambda_{i} `(\eta_{i}^{2})^{(j)}  \right\} \frac{1}{\sqrt{(\eta_{2}^{2})^{(j)} \left(r_{0}^{(j)}- (\eta_{2}^{2})^{(j)}\right)}}}{ \int_{0}^{r_{0}^{(j)}} \exp\left\{ \frac{n\beta}{2} \left( J + \frac{1}{J} \right) r_{0}^{(j)}- \frac{c\beta}{2} (\eta_{2}^{2})^{(j)} + \frac{n\beta}{2} \sum_{i=3}^{n} \lambda_{i} (\eta_{i}^{2})^{(j)}  \right\} \frac{1}{\sqrt{(\eta_{2}^{2})^{(j)} \left(r_{0}^{(j)}- (\eta_{2}^{2})^{(j)}\right)}} d(\eta_{2}^{2})^{(j)}}\\
&\propto \exp \left\{ - \frac{c\beta}{2} (\eta_{2}^{2})^{(j)}  \right\} \frac{1}{\sqrt{(\eta_{2}^{2})^{(j)} \left( r_{0}^{(j)}- (\eta_{2}^{2})^{(j)} \right)}}.
\end{split}
\end{equation}
This comes from the fact that under the uniform measure on the surface of the unit sphere, $(\eta_{1}^{2},\ldots, \eta_{n}^{2}) \sim \mathrm{Dir}\left( \frac{1}{2}, \ldots, \frac{1}{2} \right).$
Now 
\begin{equation}
    \begin{split}
     &\mu_{n,\beta}^{\mathrm{Gibbs}}\left(\mathrm{Ov}(\underline{\sigma}^{(1)},\underline{\sigma}^{(2)}) \in B \cap \underline{\sigma}^{(1)}  \in A \cap \underline{\sigma}^{(2)} \in A\right)\\
     &= \mathrm{E}_{\mu_{n,\beta}^{\mathrm{Gibbs}}}\left[ \mathbb{I}_{\mathrm{Ov}(\underline{\sigma}^{(1)},\underline{\sigma}^{(2)}) \in B} \times \mathbb{I}_{\underline{\sigma}^{(1)}  \in A \cap \underline{\sigma}^{(2)} \in A} \right]\\
     &= \mathrm{E}_{\mu_{n,\beta}^{\mathrm{Gibbs}}}\left[  \mathrm{E}\left[ \mathbb{I}_{\mathrm{Ov}(\underline{\sigma}^{(1)},\underline{\sigma}^{(2)}) \in B} \times \mathbb{I}_{\underline{\sigma}^{(1)}  \in A \cap \underline{\sigma}^{(2)} \in A}  \left|\left(r_{0}^{(j)},(\eta_{3}^{2})^{(j)},\ldots, (\eta_{n}^{2})^{(j)}\right)  \right. \right]\right]\\
     &= \mathrm{E}_{\mu_{n,\beta}^{\mathrm{Gibbs}}}\left[\mathrm{E}\left[\mathrm{E} \left[\mathbb{I}_{\mathrm{Ov}(\underline{\sigma}^{(1)},\underline{\sigma}^{(2)}) \in B} \times \mathbb{I}_{\underline{\sigma}^{(1)}  \in A \cap \underline{\sigma}^{(2)} \in A}\left|\left( (\eta_{1}^{2})^{(j)},\ldots, (\eta_{n}^{2})^{(j)} \right)  \right.\right]\left| \left(r_{0}^{(j)},(\eta_{3}^{2})^{(j)},\ldots, (\eta_{n}^{2})^{(j)}\right) \right.\right]\right]
    \end{split}
\end{equation}
Firstly, observe that 
\[
\underline{\sigma}^{(j)}\left| \left( (\eta_{1}^{2})^{(j)},\ldots, (\eta_{n}^{2})^{(j)} \right)  \right. \stackrel{d}{=} \sum_{i=1}^{n} \zeta_{i}^{(j)} \left| \eta_{i}^{(j)} \right| u_{i}
\] 
under $\mu_{n,\beta}^{\mathrm{Gibbs}}$. Here $\zeta_{i}^{(j)}$'s are i.i.d. Rademacher random variables for $1\le i \le n$ and $1\le j \le 2$. So, 
\[
\mathrm{Ov}\left( \underline{\sigma}^{(1)}, \underline{\sigma}^{(2)} \right)\left| \left( (\eta_{1}^{2})^{(j)},\ldots, (\eta_{n}^{2})^{(j)} \right) \right. \stackrel{d}{=} \sum_{i=1}^{n} \zeta_{i}^{(1)} \zeta_{i}^{(2)} \left| \eta_{i}^{(1)} \eta_{i}^{(2)} \right|
\]
 First bound the term $\sum_{i=3}^{n} \zeta_{i}^{(1)} \zeta_{i}^{(2)} \left| \eta_{i}^{(1)} \eta_{i}^{(2)}\right|$. Observe that inside $A$, $0\le r_{i} \le \tau\tilde{\varepsilon} $ for $ 1 \le i \le K$ for some known $\tau \in [0,\infty)$. This will make $(\eta_{i}^{2})^{(j)} \le \tau \tilde{\varepsilon}$ for $i \ge 3$. It is easy to see that $\mathrm{E}\left[ \sum_{i=3}^{n} \zeta_{i}^{(1)} \zeta_{i}^{(2)} \left| \eta_{i}^{(1)} \eta_{i}^{(2)}\right|\right]=0$ 
 \[
 \mathrm{Var}\left[ \sum_{i=3}^{n} \zeta_{i}^{(1)} \zeta_{i}^{(2)} \left| \eta_{i}^{(1)} \eta_{i}^{(2)}\right| \right]= \sum_{i=3}^{n} (\eta_{i}^{2})^{(1)} (\eta_{i}^{2})^{(2)} \le \tau \tilde{\varepsilon} \left(\frac{1}{\beta J}+ \tilde{\varepsilon}\right).
 \]
 So, this quantity is very small. An application of Chebyshev's inequality will prove that $\sum_{i=3}^{n} \zeta_{i}^{(1)} \zeta_{i}^{(2)} \left| \eta_{i}^{(1)} \eta_{i}^{(2)}\right|$ is very small with high probability. Finally,  
 \[
\sum_{i=1}^{2} \zeta_{i}^{(1)} \zeta^{(2)} \left| \eta_{i}^{(1)} \eta_{i}^{(2)}\right| \left| (r_{0},(\eta_{3}^{2})^{(j)},\ldots, (\eta_{n}^{2})^{(j)} ) \right.
 \] 
 is distributed as 
 \begin{equation}\label{eq:rv}
 \sqrt{Y^{(1)}Y^{(2)}} \zeta_{1} + \sqrt{(r_{0}- Y^{(1)})(r_{0}- Y^{(2)})} \zeta_{2}
 \end{equation}
 where $Y^{(1)}$ and $Y^{(2)}$ are i.i.d. with density given in \eqref{eq:densityfinale}. First take $n \to \infty$ and then take $\tilde{\varepsilon} \to 0$. Observe that the random variable in \eqref{eq:rv} converges in distribution to the random variable in the r.h.s. of \eqref{eq:resultden} as $\tilde{\varepsilon} \to 0$. This gives us the result.
\end{proof}
\begin{appendices}
\section{Proof of Lemma \eqref{lem:open}}
\begin{proof}
From \eqref{eq:maxim1} we have 
\begin{eqnarray}\label{maxim2}
F(r_{0},\ldots, r_{K})&=& \frac{n\beta}{2}\sum_{i=0}^{K}\tilde{\lambda}_ir_i+\frac{n}{2}\text{log}K+\frac{n}{2K}\sum_{i=1}^K\text{log}r_i\\
&=&\frac{n}{2}\left[\beta (J+\frac{1}{J})r_0+\beta\sum_{i=1}^{K}\tilde{\lambda}_ir_i+\frac{1}{K}\sum_{i=1}^K\text{log}r_i\right]\nonumber
\end{eqnarray}
Let $(r_{0},\ldots, r_{K})\in S_{K} \backslash S(\hat{r}_{0},\ldots , \hat{r}_{K})\text{ and }r_i=\hat{r}_i+e_{i},\,0\le i\le K$. Then \eqref{maxim2} becomes:
\begin{eqnarray}\label{eq:fvalue}
& &\frac{n}{2}\left[\beta (J+\frac{1}{J})(\hat{r}_0+e_0)+\beta\sum_{i=1}^{K}\tilde{\lambda}_i(\hat{r}_i+e_{i})+\frac{1}{K}\sum_{i=1}^K\text{log}(\hat{r}_i+e_{i})\right]\hspace{0.3cm}\nonumber\\
&=&F(\hat{r}_{0},\ldots, \hat{r}_{K})+\frac{n}{2}\left[\beta (J+\frac{1}{J})e_0+\beta\sum_{i=1}^{K}\tilde{\lambda}_ie_{i}+\frac{1}{K}\sum_{i=1}^K\text{log}(1+\frac{e_{i}}{\hat{r}_i})\right]
\end{eqnarray}
Since $(\hat{r}_{0}, \ldots, \hat{r}_{K})$ is the optimizer in Theorem \ref{thm:free}, Thus by \eqref{eq:partial} we have 
\begin{eqnarray}\label{eq:simplex}
\beta (J+\frac{1}{J})e_0+\beta\sum_{i=1}^{K}\tilde{\lambda}_ie_{i}+\frac{1}{K}\sum_{i=1}^K\frac{e_{i}}{\hat{r}_i}=0
\end{eqnarray}
Thus by \eqref{eq:simplex}, \eqref{eq:fvalue} becomes 
\begin{eqnarray}\label{eq:fvalue1}
F(\hat{r}_{0},\ldots, \hat{r}_{K})+\frac{n}{2}\left[\frac{1}{K}\sum_{i=1}^K\text{log}(1+\frac{e_{i}}{\hat{r}_i})-\frac{1}{K}\sum_{i=1}^K\frac{e_{i}}{\hat{r}_i})\right]
\end{eqnarray}
As $(r_{0},\ldots, r_{K})\in S_{K} \backslash S(\hat{r}_{0},\ldots , \hat{r}_{K})$, then $|r_{i_0}-\hat{r}_{i_0}|>\tilde{\varepsilon}$ for some $i_0\in\{0,1,\ldots,K\}$.
We split above into two cases $i_0\ne0$, $i_0=0$.
\smallskip

\textbf{Case I:} We consider $i_0\ne0$. 
Thus by equation \eqref{eq:partial}, we have $\frac{e_{i_0}}{\hat{r}_{i_0}}> \tilde{\varepsilon}K\beta ((J+\frac{1}{J})-\tilde{\lambda}_{i_0})\ge K^{1-\frac{1}{1000}}\beta ((J+\frac{1}{J})-\tilde{\lambda}_{i_0})$ and $\text{log}(1+\frac{e_{i_0}}{\hat{r}_{i_0}})= \frac{e_{i_0}}{\hat{r}_{i_0}}o(1)$, as $K$ is very large. Then we have the following
\begin{eqnarray}\label{eq:aproxlog}
\frac{1}{K}\sum_{i=1}^K\text{log}(1+\frac{e_{i}}{\hat{r}_i})-\frac{1}{K}\sum_{i=1}^K\frac{e_{i}}{\hat{r}_i} \le -K^{-\frac{1}{1000}}\beta ((J+\frac{1}{J})-\tilde{\lambda}_i)
\end{eqnarray}
From \eqref{eq:fvalue1} and \eqref{eq:aproxlog},
 we have 
\begin{eqnarray}\label{eq:case1}
F(r_{0},\ldots, r_{K})&\le& F(\hat{r}_{0},\ldots, \hat{r}_{K})- \frac{n}{2}\beta ((J+\frac{1}{J})-\tilde{\lambda}_1)\tilde{\varepsilon}\\
&\le& F(\hat{r}_{0},\ldots, \hat{r}_{K})- \frac{n}{2}\beta ((J+\frac{1}{J})-\tilde{\lambda}_1)\tilde{\varepsilon}^2 \nonumber
\end{eqnarray}
\smallskip

\textbf{Case II:} We consider $i_0=0$. Let us fix $r_0=\alpha$ and define 
\begin{eqnarray}\label{eq:maxim2}
\tilde{F}(\alpha)=\max_{(r_0=\alpha, r_1,\ldots,r_K\in S_K)} F(\alpha,r_1, \ldots, r_{K})
\end{eqnarray}
,where \[F(\alpha,r_1, \ldots, r_{K}):=\frac{n\beta}{2}\left(J+\frac{1}{J}\right)\alpha+\frac{n\beta}{2}\sum_{i=1}^{K}\tilde{\lambda}_ir_i+\frac{n}{2}\text{log}K+\frac{n}{2K}\sum_{i=1}^K\text{log}r_i
\]
One can observe that $F(\alpha,r_1, \ldots, r_{K})\le \tilde{F}(\alpha)$ for any $(\alpha,r_1, \ldots, r_{K})\in S_K$. Thus by equation \eqref{eq:maxim1} we have 
\begin{eqnarray}\label{}
F(\alpha,r_1, \ldots, r_{K})\le \tilde{F}(\hat{r}_{0})= F(\hat{r}_{0},\ldots, \hat{r}_{K})
\end{eqnarray}
\textbf{Behavior of graph of $\tilde{F}$:}
From \eqref{eq:maxim2}, putting $v_i=\frac{r_i}{1-r_0},\,1\le i\le K$, notice that the following function 
\begin{equation}
\begin{split}
&\frac{n\beta}{2}\sum_{i=1}^{K}\tilde{\lambda}_iv_i(1-r_0)+\frac{n}{2}\text{log}K+\frac{n}{2K}\sum_{i=1}^K\text{log}v_i(1-r_0)\\
&= \frac{n}{2} \log(1-r_{0})+ \frac{n\beta(1-r_{0})}{2} \sum_{i=1}^{K}\tilde{\lambda}_iv_i + \frac{n}{2} \log K  + \frac{n}{2K} \sum_{i=1}^{K} \log v_{i}
\end{split}
\end{equation}
is concave in $v_1,v_2,\ldots, v_K$. Also, using Lagrange multiplier's method and from case-I, II in \cite{Ban25} one can observe that, $\beta(1-r_0)<1$ and $\beta(1-r_0)>1$ correspond the high and low temperature case respectively.
Firstly consider $r_{0}$ such that $\beta(1-r_{0})<1$. In this case,
from equation \eqref{eq: optvalue}, we have 
\begin{eqnarray}\label{eq:tildef}
\tilde{F}(r_0)&=&\frac{n\beta}{2}\left(J+\frac{1}{J}\right)r_{0}+\frac{n}{2}\text{log}(1-r_0)+\frac{n\beta^2(1-r_0)^2}{4}\nonumber\\ \frac{d\tilde{F}}{dr_0}&=&\frac{n\beta}{2}\left(J+\frac{1}{J}\right)-\frac{n}{2}\frac{1}{1-r_0}-\frac{n\beta^2(1-r_0)}{2} \nonumber\\
\frac{d^2\tilde{F}}{dr^2_0}&=&-\frac{n}{2}\frac{1}{(1-r_0)^2}+\frac{n\beta^2}{2}\hspace{0.5cm}
\end{eqnarray}
Above shows that $\tilde{F}$ has local maximum at $r_0=1-\frac{1}{\beta J}$.
\smallskip

 For the low temperature case (i.e. $r_0<1-\frac{1}{\beta}$) from equation $(2)$ in \cite{Ban25} and \eqref{eq:maxim2} we have 
\[\tilde{F}(r_0)=\frac{n\beta}{2}\left[\left(J+\frac{1}{J}\right)-2\right]r_0+n\beta-\frac{3n}{4}-\frac{n}{2}\text{log}(\beta)\] is an affine linear and  increasing function as $\frac{n\beta}{2}\left[\left(J+\frac{1}{J}\right)-2\right]>0.$ Thus, $\tilde{F}$ has global maximum at $r_0=1-\frac{1}{\beta J}=\hat{r}_0$
\smallskip

Now, for $(r_{0},\ldots, r_{K})\in S_{K} \backslash S(\hat{r}_{0},\ldots , \hat{r}_{K})$ with $|r_0-\hat{r}_0|>\tilde{\varepsilon}$ and from graph of $\tilde{F}$ we have the following
\begin{eqnarray}\label{eq:globalmax}
F(r_0,r_1, \ldots, r_{K})\le \text{max}\{ \tilde{F}(\hat{r}_{0}+\tilde{\varepsilon}),\,\tilde{F}(\hat{r}_{0}-\tilde{\varepsilon})\}.
\end{eqnarray}
We take $\tilde{\varepsilon}$ so that $\beta(1 - \hat{r}_{0} \pm \tilde{\varepsilon})<1$.
Now, by equation \eqref{eq: optvalue} we have 
\begin{eqnarray}\label{eq:difference}
\tilde{F}(\hat{r}_0\pm \tilde{\varepsilon})&=&\frac{n\beta}{2}\left(J+\frac{1}{J}\right)(\hat{r}_0\pm \tilde{\varepsilon})+\frac{n}{2}\text{log}(1-\hat{r}_0\mp \tilde{\varepsilon})+\frac{n\beta^2(1-\hat{r}_0\mp \tilde{\varepsilon})^2}{4}\nonumber\\
&=&\tilde{F}(\hat{r}_0)+\frac{n\beta}{2}\left(J+\frac{1}{J}\right)\tilde{\varepsilon}+\frac{n}{2}\text{log}(1-\hat{r}_0 \mp\tilde{\varepsilon})-\frac{n}{2}\text{log}(1-\hat{r}_0)\nonumber\\
&+&\frac{n\beta^2(1-\hat{r}_0\mp \tilde{\varepsilon})^2}{4}-\frac{n\beta^2(1-\hat{r}_0)^2}{4}\nonumber\\
&=&\tilde{F}(\hat{r}_0)+\frac{n\beta}{2}\left(J+\frac{1}{J}\right)\tilde{\varepsilon}+\frac{n}{2}\text{log}(1\mp\frac{\tilde{\varepsilon}}{1-\hat{r}_0})\\
&\mp&\frac{n\beta^2\tilde{\varepsilon}(2(1-\hat{r}_0)
\mp\tilde{\varepsilon})}{4}\nonumber
\end{eqnarray}
We know that $\log(1+x)-x = -\frac{x^2}{2} +O(x^3)$. A careful analysis shows that in \eqref{eq:difference} the linear term of $\tilde{\varepsilon}$ cancels out and the coefficient of the quadratic term is negative. Thus, by equation \eqref{eq:difference} we have
\begin{eqnarray}\label{eq:case2}
\tilde{F}(\hat{r}_0\pm\tilde{\varepsilon})&\le&\tilde{F}(\hat{r}_0)- \gamma n \tilde{\varepsilon}^2
\end{eqnarray}
for some fixed $\gamma \in (0,\infty)$.
This concludes the proof.

\end{proof}
\end{appendices}

\noindent
\textbf{Acknowledgements:} We are grateful to Prof. Soumendu Sundar Mukherjee and Prof. Riddhipratim Basu for private communications.
\bibliographystyle{alpha}
\bibliography{main}
\end{document}